\documentclass[openany]{book}

\ifx\directlua\undefined\ifx\XeTeXcharclass\undefined
  \else\RequirePackage[no-math]{fontspec}[2017/03/31]\fi 
  \else\RequirePackage[no-math]{fontspec}[2017/03/31]\fi 
\usepackage[medium]{dgruyter}
\usepackage{subcaption}
\usepackage{wrapfig}
\usepackage{cleveref}
\usepackage{todonotes}

\DeclareMathOperator*{\argmin}{arg\,min}

\contribution

\contributionauthor[1]{Alexander Heinlein}
\contributionauthor[2]{Amanda A. Howard}
\contributionauthor[3]{Damien Beecroft}
\contributionauthor[2]{Panos Stinis}
\affil[1]{Delft University of Technology, Delft Institute of Applied Mathematics, 2628 CD Delft, Netherlands}
\affil[2]{Pacific Northwest National Laboratory, P.O. Box 999, Richland, WA, 99352, USA}
\affil[3]{University of Washington, Seattle, WA, 98195, USA}
\runningauthor{Heinlein, Howard, Beecroft, Stinis}
\contributiontitle{Multifidelity domain decomposition-based physics-informed neural networks and operators for time-dependent problems}
\runningtitle{Multifidelity DD NNs}
\abstract{
	Multiscale problems are challenging for neural network-based discretizations of differential equations, such as physics-informed neural networks (PINNs). This can be (partly) attributed to the so-called spectral bias of neural networks. To improve the performance of PINNs for time-dependent problems, a combination of multifidelity stacking PINNs and domain decomposition-based finite basis PINNs is employed. In particular, to learn the high-fidelity part of the multifidelity model, a domain decomposition in time is employed. The performance is investigated for a pendulum and a two-frequency problem as well as the Allen-Cahn equation. It can be observed that the domain decomposition approach clearly improves the PINN and stacking PINN approaches. Finally, it is demonstrated that the FBPINN approach can be extended to multifidelity physics-informed deep operator networks.
}

\begin{document}

\makecontributiontitle

\section{Introduction}

Many problems arising in science and engineering exhibit a multiscale nature, with different processes taking place on various temporal and spatial scales. The solution of these problems is generally difficult for numerical methods. Multiscale methods have been developed to make numerical simulations of such problems feasible; examples are the multiscale finite element~\cite{hou_multiscale_2009}, homogeneous multiscale~\cite{e_heterogeneous_2007}, generalized finite element~\cite{babuska_generalized_1983} or variational multiscale~\cite{hughes_multiscale_1995} methods. 

In recent years, inspired by the early work by Lagaris et al.~\cite{lagaris_artificial_1998}, machine learning based techniques for the solution of partial differential equations (PDEs) have been developed. In this paper, we focus on physics-informed neural networks (PINNs)~\cite{raissi_physics-informed_2019}; other methods, such as the Deep Ritz method~\cite{e_deep_2018}, have been developed around the same time. Those methods have many potential advantages: they are easy to implement, allow for direct integration of data, and can be employed to solve inverse and high-dimensional problems. However, their convergence properties are not yet well-understood, and hence, the resulting accuracy is often limited. This can be partly accounted to the spectral bias of neural networks (NNs)~\cite{rahaman_spectral_2019}, meaning that NNs tend to learn low frequency components of functions much faster than high frequency components. In multiscale problems, the high frequency components typically correspond to the fine scales, whereas the low frequency components correspond to the coarse scales. Therefore, multiscale problems are also particularly challenging to solve using PINNs.

In this paper, we aim at combining two techniques that have recently been developed to improve the training of PINNs in this context. On the one hand, we consider the multi-fidelity training approach introduced for PINNs~\cite{meng2020composite} and extended to Deep Operator Networks (DeepONets,~\cite{lu_deeponet_2021}) in~\cite{howard_multifidelity_2023, lu2022multifidelity, de2022bi}. In particular, we consider the approach of stacked PINNs~\cite{howard_stacked_2023} in which multiple networks are stacked on top of each other, such that models on top of the stack may learn those features that are not captured by the previous models. On the other hand, we employ multi-level Schwarz domain decomposition neural network architectures~\cite{dolean_multilevel_2023}, which are based on the finite-basis PINNs (FBPINNs)~\cite{moseley_finite_2023} approach. In this approach, the learning of multiscale features is improved by localization. In particular, the network architecture is decomposed, such that the individual parts of the network learn features on the corresponding spatial or temporal scale. For an overview on the combination of domain decomposition approaches and machine learning, see, for instance,~\cite{heinlein_combining_2021}.

In related recent works, methods for iteratively training PINNs to progressively reduce the errors have been developed; see~\cite{ainsworth2021galerkin,ainsworth2022galerkin,aldirany2023multi,wang2023multi}. These approaches, which vary in their implementation details, train each new network to reduce the residual from the previous network. In contrast, the work presented here trains for the entire solution at each iteration. 

This paper is structured as follows: first, in~\Cref{section:method}, we describe the methodological framework. In particular, we first discuss PINNs in~\Cref{section:pinns}, then we describe multifidelity stacking PINNs in~\Cref{section:mfpinns}, as well as the domain decomposition approach in~\Cref{section:dd}. Next, we introduce the specific domain decomposition in time employed in the model problems in~\Cref{section:dd_time}. In~\Cref{section:results}, we present numerical results for several model problems, a pendulum and a two-frequency problem as well as the Allen-Cahn equation. We extend the results to DeepONets in~\Cref{section:results-deeponets}. We conclude with a brief discussion of the current and future work in \Cref{section:discussion}. All training parameters used to generate the results are given in~\Cref{tab:parameters_PINNs}.


\section{Methodology} \label{section:method}

\subsection{Physics-informed neural networks} \label{section:pinns}

We consider a generic differential equation-based problem in residual form: find $u$ such that
\begin{equation} \label{eq:problem}
	\begin{aligned}
		\mathcal{A} u & = 0 \quad \text{in } \Omega, \\
		\mathcal{B} u & = 0 \quad \text{on } \partial\Omega, \\
	\end{aligned}
\end{equation}
where $\mathcal{A}$ is a differential operator and $\mathcal{B}$ an operator for specifying the initial or boundary conditions. The solution $u$ is defined on the domain $\Omega$ and should have sufficient regularity to apply $\mathcal{A}$ and $\mathcal{B}$. 
In order to solve~\cref{eq:problem}, we follow~\cite{lagaris_artificial_1998} and employ a collocation approach. In particular, we exploit that solving~\cref{eq:problem} is equivalent to solving
$
	\argmin_{\mathcal{B} u = 0 \text{ on } \partial\Omega} 
	\,
	\int_\Omega \left( \mathcal{A} u (\mathbf{x}) \right)^2 \, d\mathbf{x}.
$
We discretize the solution using a neural network $\hat u (\mathbf{x}, \theta)$, with parameters $\theta$, and the integral is approximated by the sum
$$
	\argmin_{\mathcal{B} \hat u (\mathbf{x}, \theta) = 0 \text{ on } \partial\Omega}
	\, 
	\sum_{\mathbf{x}_i \in \Omega} \left( \mathcal{A} \hat u (\mathbf{x}_i, \theta) \right)^2,
$$
where the collocation points $\mathbf{x}_i$ are sampled from $\Omega$.
Different types of neural network architectures may be employed, and we will employ a combination of the approaches explained in~\cref*{section:mfpinns,section:dd}.

The initial or boundary conditions in the second equation of~\cref{eq:problem} can be enforced via hard or soft constraints. In the approach of hard constraints, they are explicitly implemented in the neural network function; cf.~\cite{lagaris_artificial_1998}. In this paper, we employ the approach of soft constraints instead, in which we incorporate the constraints into the loss function:
\begin{equation} \label{eq:loss}
	\argmin_{\theta}
	\,
	\lambda_r
	\sum_{\mathbf{x}_i \in \Omega} \left( \mathcal{A} \hat u (\mathbf{x}_i, \theta) \right)^2
	+
	\lambda_{bc}
	\sum_{\mathbf{x}_i \in \partial\Omega} \left( \mathcal{B} \hat u (\mathbf{x}_i, \theta) \right)^2 
\end{equation}
Here, $\lambda_r$ and $\lambda_{bc}$ weight the residual of the differential equation and the initial and boundary conditions in the loss function. As discussed, for instance, in~\cite{wang_when_2022} an appropriate weighting is crucial for the convergence in optimizing~\cref{eq:loss} using a gradient-based optimization method. $\theta$ denotes all the trainable parameters in the network. 

This approach has also been denoted as physics-informed neural networks (PINNs) in~\cite{raissi_physics-informed_2019}.

\subsection{Multifidelity stacking PINNs} \label{section:mfpinns}
Multifidelity PINNs use two NNs to learn the correlation between low- and high-fidelity physics \cite{meng2020composite}. The goal is to train a linear network (with no activation function) 
to learn the linear correlation between the low- and high-fidelity models, and a nonlinear network to learn the nonlinear correlation. By training a linear network, the resulting model is more expressive than just assuming that the 
correlation between the models is the identity. Moreover, under the assumption that the main part of the correlation is linear, separating the network into the linear and nonlinear parts allows for a smaller nonlinear network. 

To train a multifidelity PINN we first train a standard single fidelity PINN $\hat u^{SF}(\mathbf{x}, \theta^{SF})$. In a second step, we then train a multifidelity network $\hat u^{MF}$, which consists of linear and nonlinear subnetworks that learn the correlation between the single fidelity PINN $\hat u^{SF}(\mathbf{x}, \theta^{SF})$ and the solution: 
\begin{equation} \label{eq:multifelity_network_standard}
	\hat u^{MF}(\mathbf{x}, \theta^{MF})=(1-\left|\alpha\right|)\hat u_{linear}^{MF}(\mathbf{x}, \hat u^{SF}, \theta^{MF})+ \left|\alpha\right| \hat u_{nonlinear}^{MF}(\mathbf{x}, \hat u^{SF}, \theta^{MF}).
\end{equation}
The linear network does not have activation functions, to force learning a linear correlation, and can be very small. 
$\alpha$ is a trainable parameter to enforce maximizing the linear correlation.

The loss function in~\cref{eq:loss} is modified to include the penalty $\alpha^4$: 
\begin{equation} \label{eq:loss_MF}
	\argmin_{\theta}
	\,
	\lambda_{r}
	\sum_{\mathbf{x}_i \in \Omega} \left( \mathcal{A} \hat u (\mathbf{x}_i, \theta) \right)^2
	+
	\lambda_{bc}
	\sum_{\mathbf{x}_i \in \partial\Omega} \left( \mathcal{B} \hat u (\mathbf{x}_i, \theta) \right)^2  + \lambda_\alpha \alpha^4
\end{equation}

In multifidelity stacking PINNs as presented in~\cite{howard_stacked_2023}, multifidelity PINNs are trained recursively, 
each taking the output of the previously trained stacking layer as input. In this way, the previous layer serves as the low-fidelity model for the new stacking layer. The difference between \cite{howard_stacked_2023} and the current work is that \cite{howard_stacked_2023} does not consider domain decomposition, so each stacking layer has a single multifidelity PINN covering the entire domain. The approach considered here is more flexible and, as we will show in Section \ref{section:results}, results in smaller relative errors when trained on the same equations.

\subsection{Domain decomposition-based neural network architectures} \label{section:dd}

It has been observed in~\cite{moseley_finite_2023} that the high frequency components in the solution can be learned better if a domain decomposition is introduced into the PINN approach. To scale to larger numbers of subdomains, this approach has first been extended to two-levels in~\cite{dolean_finite_2022} and then to an arbitrary number of levels in~\cite{dolean_multilevel_2023}.
The general idea of the domain decomposition-based finite basis PINNs (FBPINNs) is to decompose the computational domain $\Omega$ into $J$ overlapping subdomains $\Omega_j$, $\Omega = \cup_{j=1}^J \Omega_j$. As before, $\Omega$ may be a space-time domain, and in this work we will focus on domain decomposition in time. On each subdomain, we define a space of network functions 
$\mathcal{V}_j = \left\{ \hat u_j(\mathbf{x}, \theta_j) | \mathbf{x} \in \Omega_j, \theta_j \in \Theta_j \right\},$
where ${\hat u}_j(\mathbf{x}, \theta_j)$ denotes a PINN model, $\Theta_j = \mathbb{R}^{k_j}$ is the space of all trainable neural network parameters, and $k_j$ is the number of network parameters. 

In order to represent the global solution of a given problem, we define window functions $\omega_j$ with $\text{supp}(\omega_j) \subset \Omega_j$ such that $\{\omega_j\}_{j=1}^J$ form a partition of unity, that is,
$
	\sum_{j=1}^J \omega_j = 1 \text{ on } \Omega.
$
Then, we can define a global neural network space $\mathcal{V} = \sum_{j=1}^J \omega_j \mathcal{V}_j $, and the global FBPINN function reads 
$
	\hat u(\mathbf{x}, \theta) = \sum_{j=1}^J \omega_j \hat u_j(\mathbf{x}, \theta_j).
$
It has been observed that this approach may significantly improve the performance of PINNs; cf.~\cite{moseley_finite_2023}. However, similar to classical domain decomposition methods~\cite{toselli_domain_2005}, the one-level approach is not scalable to large numbers of subdomains; see~\cite{dolean_finite_2022,dolean_multilevel_2023}.

\begin{figure}[t]
	\begin{center}
		\begin{tikzpicture}
			\draw[>=stealth,|<->|] 				( 0.0, 0.0) -- node[below=-0.5mm] {\small $\Omega$} +(10.0, 0);
			
			\draw[red!30!black,very thick] 		( 0.0, 0.6) -- node[below=-0.5mm] {\small $\Omega_1^{(3)}$} +( 3.0, 0);
			\draw[red!50!black,very thick] 		( 2.0, 0.7) -- node[below=-0.5mm] {\small $\Omega_2^{(3)}$} +( 3.5, 0);
			\draw[red!70!black,very thick] 		( 4.5, 0.6) -- node[below=-0.5mm] {\small $\Omega_3^{(3)}$} +( 3.5, 0);
			\draw[red!90!black,very thick] 		( 7.0, 0.7) -- node[below=-0.5mm] {\small $\Omega_4^{(3)}$} +( 3.0, 0);
			
			\draw[blue!40!black,very thick] 	( 0.0, 1.3) -- node[below=-0.5mm] {\small $\Omega_1^{(2)}$} +( 6.0, 0);
			\draw[blue!80!black,very thick] 	( 4.0, 1.4) -- node[below=-0.5mm] {\small $\Omega_2^{(2)}$} +( 6.0, 0);
			
			\draw[green!60!black,very thick] 	( 0.0, 2.0) -- node[below=-0.5mm] {\small $\Omega_1^{(1)}$} +(10.0, 0);
			
		\end{tikzpicture}
	\end{center}
	\caption{Multilevel overlapping domain decomposition of $\Omega$ with $L=3$ levels.
		\label{fig:dd}
	}
\end{figure}
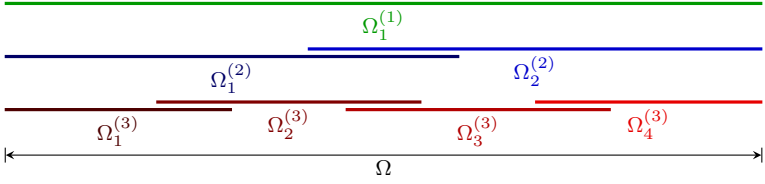
To improve the scalability and the performance for multiscale problems, a hierarchy of domain decompositions may be employed. Define $L$ levels of domain decompositions, with the overlapping domain decomposition at level $l$ denoted by $D^{(l)} = \{\Omega_j^{(l)}\}_{j=1}^{J^{(l)}}$, where $\Omega = \cup_{j=1}^{J^{(l)}} \Omega_j^{(l)}$ and $J^{(l)}$ is the number of subdomains at level $l$; cf.~\Cref{fig:dd}. Even though there is generally no restriction on the overlapping domain decompositions, we choose $J^{(1)} = 1$, so the first level corresponds to a single global subdomain, and $J^{(l)} < J^{(l+1)}$ for all $l = 1, \ldots, L$. 

Now, on each level $l$ we define window functions $\omega_j^{(l)}$ to be a partition of unity, so $\sum_{j=1}^{J^{(l)}}\omega_j^{(l)} =1$, and $\text{supp}(\omega_j^{(l)}) \subset \Omega_j^{(l)}$. Similar to the one-level case, this yields the global neural network space
$ 
	\mathcal{V} 
	= 
	\sum_{l=1}^L\sum_{j=1}^{J^{(l)}} \omega_j^{(l)} \mathcal{V}_j^{(l)}
$
and the global network function defined in terms of $\theta = \cup_{l=1}^{L} \theta^{(l)}$ and $\theta^{(l)} = \cup_{j=1}^{J^{(l)}} \theta_j^{(l)}$:
\begin{equation} 
	\hat u(\mathbf{x}, \theta) 
	= 
	\frac{1}{L}
	\sum_{l=1}^{L} \hat u^{(l)}(\mathbf{x}, \theta^{(l)})
	\quad \text{with} \quad
	\hat u^{(l)}(\mathbf{x}, \theta^{(l)}) 
	= 
	\sum_{j=1}^{J^{(l)}} \omega_j^{(l)} \hat u_j^{(l)}(\mathbf{x}, \theta_j^{(l)}). \label{eq:sum_FBPINN}
\end{equation} 
 It has been observed in~\cite{dolean_finite_2022,dolean_multilevel_2023} that,
 due to increased communication between the subdomain models, the multilevel FBPINN approach may significantly improve the performance over the one-level approach.

\subsection{Stacking FBPINNs} \label{section:stacking_fbpinns}

We combine the multifidelity stacking PINNs and the FBPINNs as follows: In the first level, we train a standard single fidelity PINN across the full domain $\Omega^{(0)} = \Omega$. Then, for each level $l > 0$ we use a FBPINN network architecture modified to consist of multifidelity networks that takes as input the network from the previous level $l-1$:
\begin{equation}
   	\hat u^{(l)}(\mathbf{x}, \theta^{(l)}) 
	= 
	\sum_{j=1}^{J^{(l)}} \omega_j^{(l)} \hat u_{j, MF}^{(l)}(\mathbf{x},  \hat u^{(l-1)}, \theta_j^{(l)}).  \label{eq:stackingFBPINN}
\end{equation}
We note that~\Cref{eq:stackingFBPINN} differs from~\Cref{eq:sum_FBPINN} by a factor of $1/L$, because the output of FBPINNs is the sum of the networks trained at all levels, while the output of stacking FBPINNs is the sum of the networks for the final level. 
The networks at level $l$ learn the correlation between the output of the $l-1$ level and the solution, and take as input the previously learned solution $\hat u^{(l-1)}(\mathbf{x}, \theta^{(l-1)})$: 
\begin{equation*}
	\hat u_{j, MF}^{(l)}(\mathbf{x}, \theta_j^{(l)})=(1-|\alpha|)\hat u_{j, lin}^{(l)}\left(\mathbf{x}, \hat u ^{(l-1)},\theta_j^{(l)}\right)+|\alpha|\hat u_{j, nonlin}^{(l)}\left(\mathbf{x}, \hat u^{(l-1)}, \theta_j^{(l)}\right)
\end{equation*}

\section{Domain decomposition in time} \label{section:dd_time}
\begin{wrapfigure}[7]{r}{0.45\textwidth}
    \centering
        \vspace{-2pt}
            \includegraphics[width=0.4\textwidth]{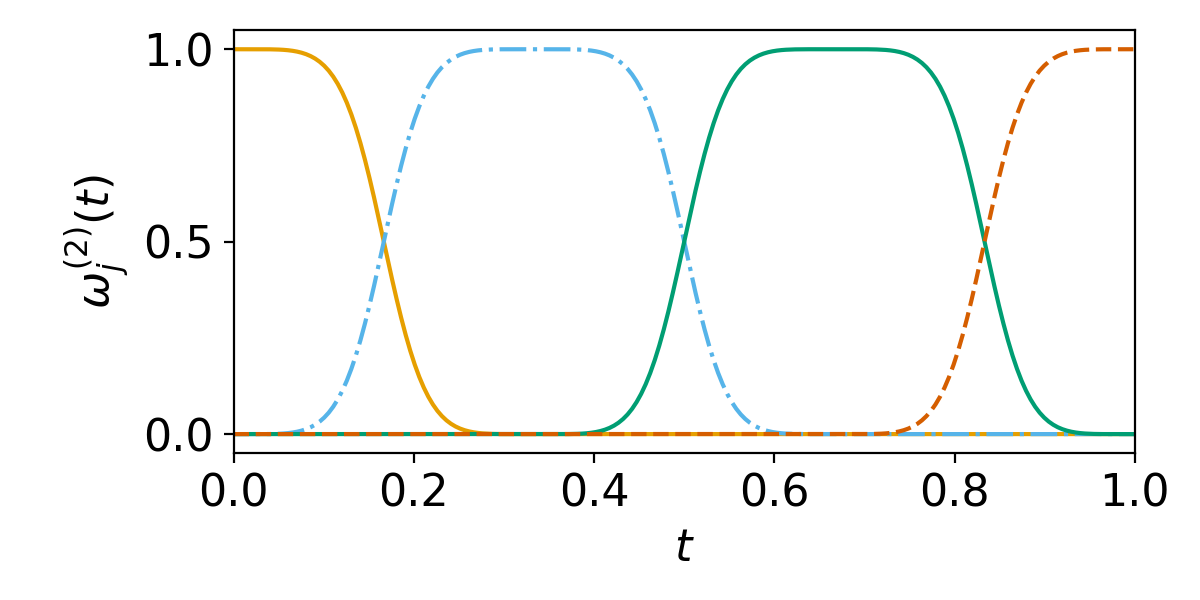}
    \caption{Window functions $\omega_j$ for $l = 2$ and $T = 1$.}
    \label{fig:weight_functions}
\end{wrapfigure}

In this work, we are particularly interested in cases were classical PINNs fail to learn the temporal evolution, 
such as a damped pendulum and the Allen-Cahn equation. We consider a domain $\Omega = \mathbf{X} \times [0, T]$ where $\mathbf{X}$ denotes the spatial domain and $T \in \mathbb{R}$. Therefore, for the stacking FBPINN approach, we consider 
the domain decomposition in time:
$$
	\Omega_j^{(l)} = \left[\frac{(j-1)T -\delta T/2}{J^{(l)}-1}, \frac{(j-1)T +\delta T/2}{J^{(l)}-1} \right]
$$
where $\delta>1$ is the overlap ratio. For $l=0$, we take 
$\Omega_1^{(0)} = [0.5T-\delta T/2, 0.5T+\delta T/2].$
The partition of unity functions are given by 
$
	\omega_j^{(l)} = \frac{\hat\omega_j^{(l)}}{\sum_{j=1}^{J^{(l)}} \hat\omega_j^{(l)}},
$
where
\begin{equation}
	\hat\omega_j^{(l)}(t) = \begin{cases}
		1 & l = 0, \\
		\left[ 1+cos\left(\pi(t-\mu_j^{(l)})/\sigma_j^{(l)}\right)\right]^2 & l>0,
	\end{cases} \label{eq:PoU}
\end{equation}
$\mu_j^{(l)} = T(j-1)/(J^{(l)}-1)$, and $\sigma_j^{(l)} = (\delta T/2)/(J^{(l)}-1)$. For simplicity, we take $J^{(l)} = 2^l$ in each case. An illustration of the window functions for $T=1$ and $l = 2$ ($J^{(2)} = 4$) is given in~\Cref{fig:weight_functions}.

As set up, each network only covers a small part of the time domain. To ease training, we scale the input in each domain to be in the range [-1, 1] by using a scaled time $\hat t = t(l-1)/T-j$ as the input to network $j$. This scaling improves the robustness of the training. 



In our applications we calculate the relative $\ell_2$ error 
$
    \frac{|| u(\mathbf{x}) - \hat u(\mathbf{x}, \theta) ||_2}{|| u(\mathbf{x})||_2}
$
where $u$ denotes the exact solution and $\hat u$ denotes the output from the multifidelity FBPINN.

\section{Results} \label{section:results}

\subsection{Pendulum} \label{sec:pendulum}

\begin{figure}[h!]
    \centering
    \begin{subfigure}[t]{0.5\textwidth}
    	\centering
    	\includegraphics[width=\textwidth]{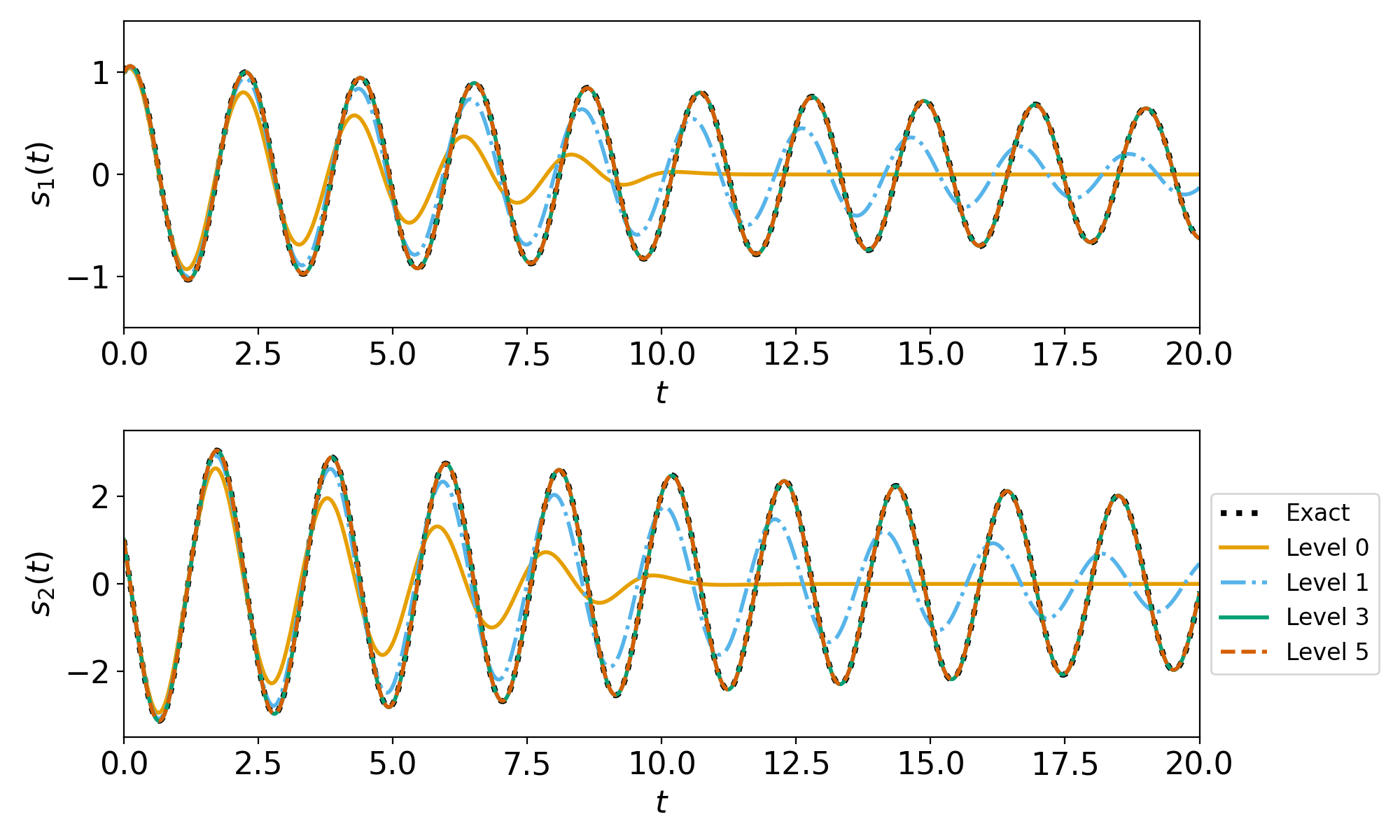}
    \end{subfigure}
    \begin{subfigure}[t]{0.35\textwidth}    
    	\centering
        \includegraphics[width=\textwidth]{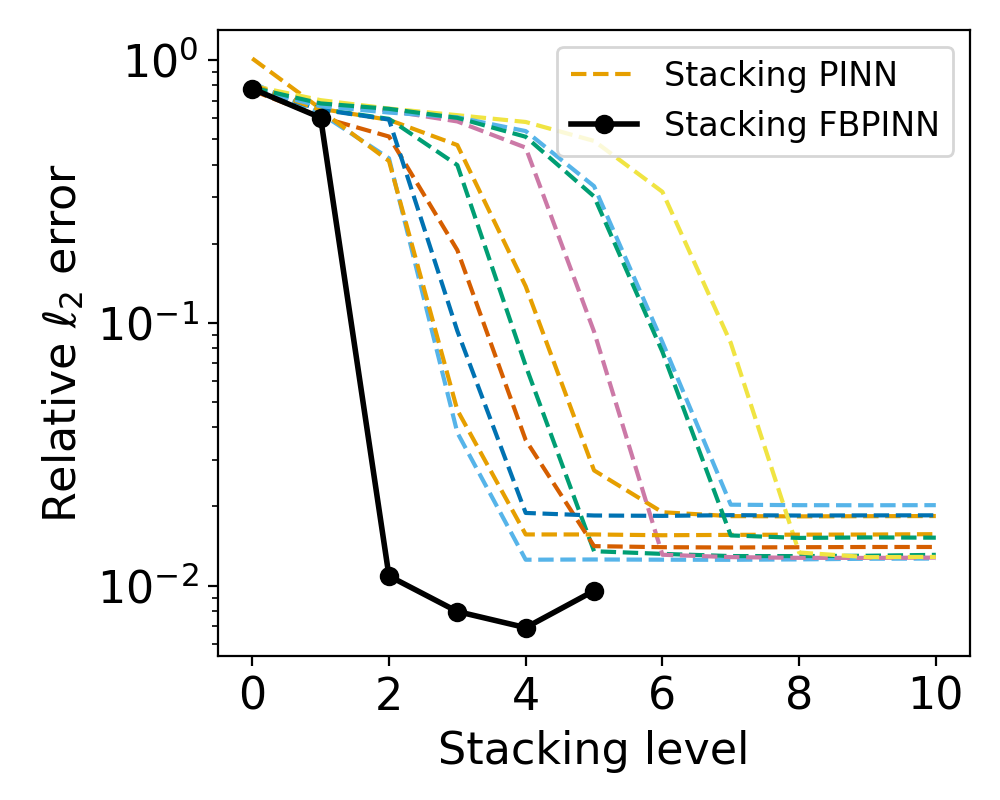}
    \end{subfigure}
    \caption{Stacking FBPINN results for the pendulum problem: \textbf{Left:} Stacking FBPINN results for an illustrative example of $s_1$ (top) and $s_2$ (bottom) as a function of time for the pendulum problem up to five stacking FBPINN levels. \textbf{Right:} Pendulum relative $\ell_2$ training errors comparing the work in the current paper (solid line) with the approach from~\cite{howard_stacked_2023} (dashed lines).}
    \label{fig:pendulum_stacking}
\end{figure}

While a relatively simple system, accurately training a PINN to predict the movement of a pendulum for long times presents challenges \cite{wang2023long}. The pendulum movement is governed by a system of two first-order ODEs for $t \in [0, T]$ 
\begin{align}
    \frac{d s_1}{dt} &= s_2, \label{eq:pendulum_1}\\
    \frac{d s_2}{dt} &= -\frac{b}{m} s_2 - \frac{g}{L} \sin(s_1), \label{eq:pendulum_2}
\end{align}
where $s_1$ and $s_2$ are the position and velocity of the pendulum, respectively. We employ the same parameters used in~\cite{wang2023long}, that is, $m=L=1$, $b=0.05$, $g=9.81$, and $T=20$. We take $s_1(0) = s_2(0) = 1$. We compare the results with those for the stacking PINN from~\cite{howard_stacked_2023}, which uses the same multifidelity architecture but only a single PINN on each level. 
As shown in~\Cref{fig:pendulum_stacking}, the stacking FBPINN is able to reach a significantly lower relative $\ell_2$ error. In addition, each network in the stacking FBPINN is significantly smaller than the networks used in the stacking PINN with the result that, at three stacking layers, the stacking FBPINN reaches a relative $\ell_2$ error of $7.4{\cdot}10^{-3}$ with only 34\,570 trainable parameters. In comparison, the best case stacking PINN from~\cite{howard_stacked_2023} requires four stacking levels to reach a relative $\ell_2$ error of $1.3{\cdot}10^{-2}$ with 63\,018 trainable parameters. 

\subsection{Multiscale problem}\label{sec:Multiscale}
We now consider a toy model problem with a low and high frequency component, inspired by~\cite{moseley_finite_2023}:
\begin{align*}
    \frac{ds}{dx} 
    & = 
    \omega_1\cos (\omega_1 x) + \omega_2\cos(\omega_2 x),\\
    s(0) &= 0,
\end{align*}
on domain $\Omega = [0, 20]$ with $\omega_1 = 1$ and $\omega_2 = 15$. The exact solution for this problem is $s(x) = \sin(\omega_1 x) + \sin(\omega_2 x)$.

\begin{figure}[h]
    \centering
    \begin{subfigure}[t]{0.5\textwidth}
    	\centering
        \includegraphics[width=\textwidth]{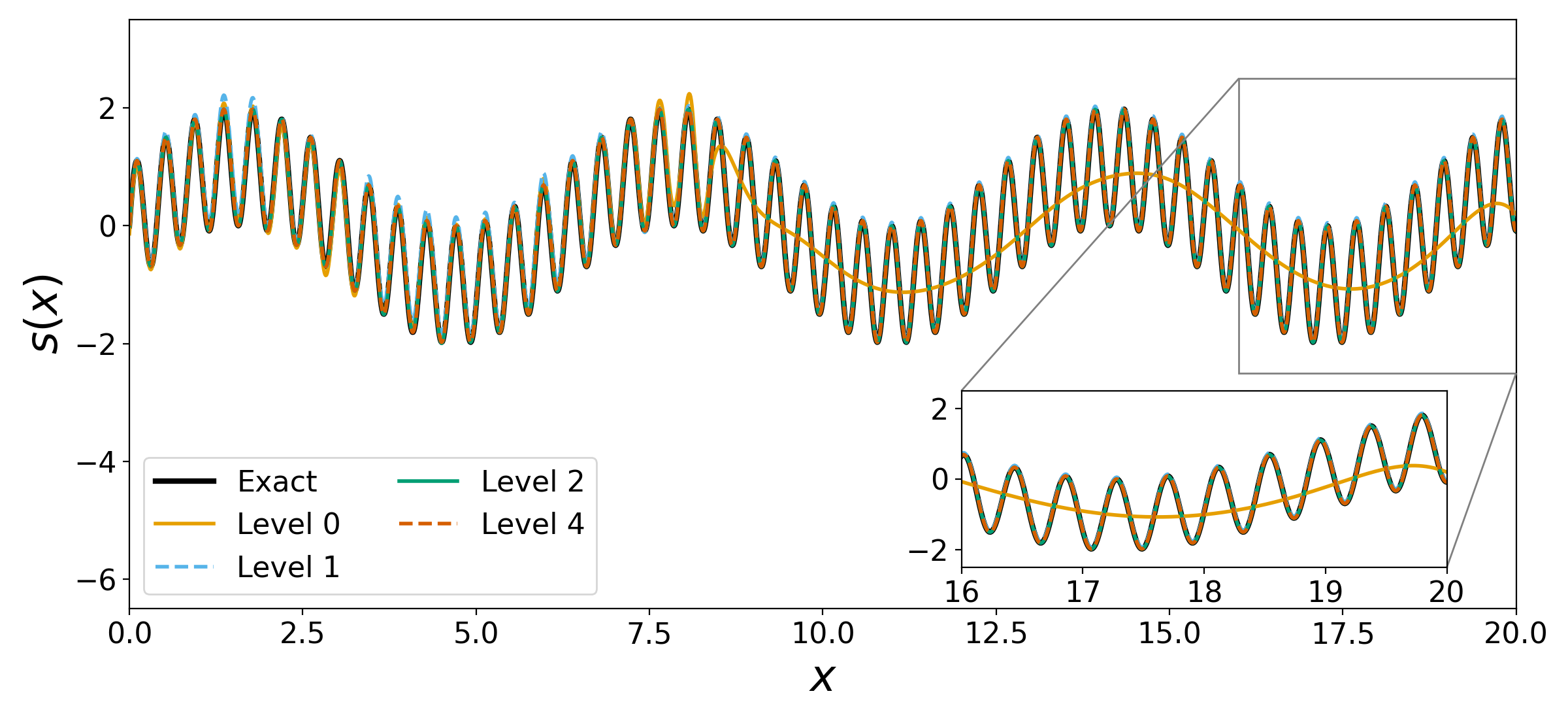}
    \end{subfigure}
    \begin{subfigure}[t]{0.35\textwidth}    
    	\centering
        \includegraphics[width=0.8\textwidth]{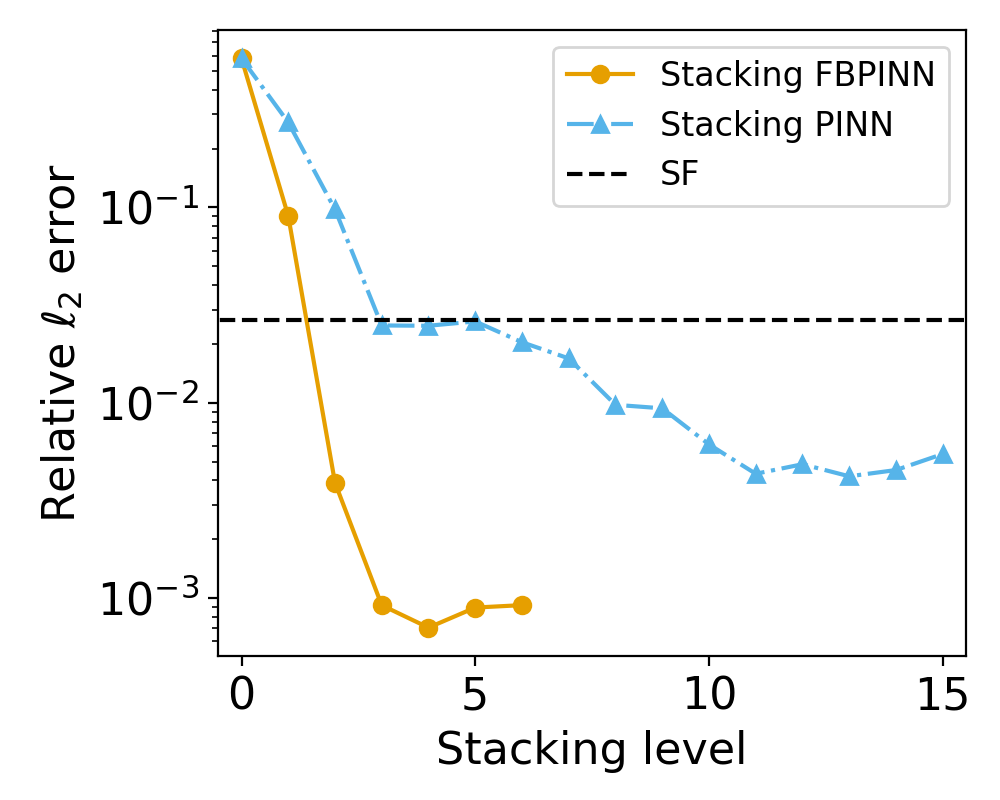}
    \end{subfigure}
    \caption{Stacking FBPINN results for the multiscale problem: \textbf{Left:} Stacking FBPINN results for the single fidelity level 0 and the first four stacking FBPINN levels. \textbf{Right:} Multiscale relative $\ell_2$ training errors comparing the work in the current paper with~\cite{howard_stacked_2023}.}
    \label{fig:multiscale_stacking}
\end{figure}

The results are shown in~\Cref{fig:multiscale_stacking}. After two stacking levels, the stacking FBPINN reaches a relative $\ell_2$ error of $4.2{\cdot}10^{-3}$,
with 7822 trainable parameters. A comparable relative $\ell_2$ error of $6.1{\cdot}10^{-3}$ is reached after 10 stacking levels with a stacking PINN with 11\,179 trainable parameters. Also shown in~\Cref{fig:multiscale_stacking} (right) is the best case SF network from \cite{howard_stacked_2023}, which has a relative $\ell_2$ error of $9.5{\cdot}10^{-2}$ with 16\,833 trainable parameters. The stacking FBPINN outperforms the SF PINN with an error more than an order of magnitude lower, with less than half the trainable parameters. 
Additionally, the final stacking FBPINN reaches a relative $\ell_2$ error of $8.3{\cdot}10^{-4}$, an order of magnitude lower than the final stacking PINN.  

\subsection{Allen-Cahn equation}\label{sec:AC}
Our third example is based on the Allen-Cahn equation and is given by 
\begin{align*}
    s_t - 0.0001 s_{xx} +5s^3-5s & = 0, & & t\in(0, 1], x \in[-1, 1], \\
    s(x, 0) & = x^2 \cos(\pi x), & & x\in[-1, 1], \\
    s(x, t) & = s(-x, t), & & t \in[0, 1], x=-1, x=1, \\
    s_x(x, t) & = s_x(-x, t), & & t\in[0, 1], x=-1, x=1. 
\end{align*}
The Allen-Cahn equation presents difficulties for PINNs when attempting to learn the full solution from $t = 0$ to $1$ with a single PINN; see, e.g.,~\cite{wight2020solving, mattey2022novel, rohrhofer2022role}.

\begin{figure}[h!]
    \centering
    \begin{subfigure}[t]{0.65\textwidth}
    	\centering
        \includegraphics[width=\textwidth]{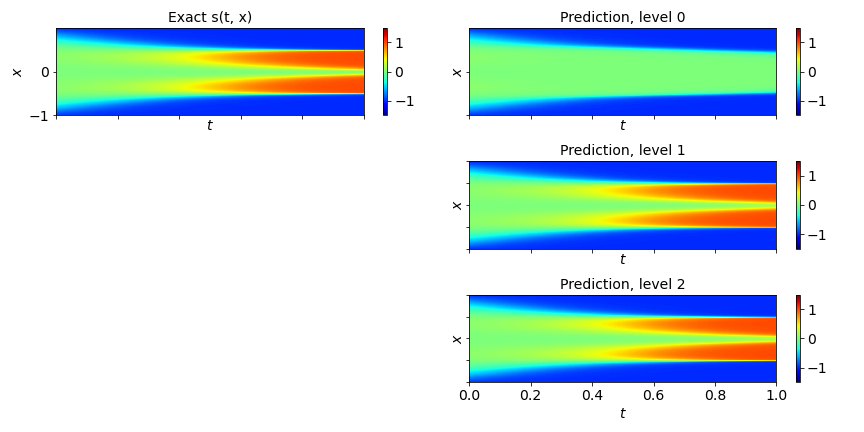}
    \end{subfigure}
    \begin{subfigure}[t]{0.2\textwidth}    
    	\centering
        \includegraphics[width=\textwidth]{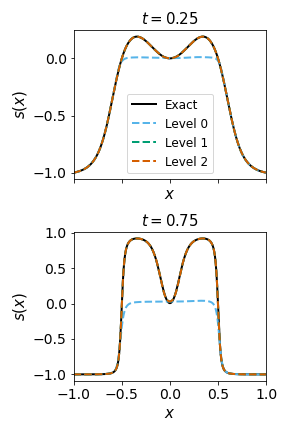}
    \end{subfigure}
    \caption{Stacking FBPINN results for the Allen-Cahn equation. \textbf{Left:} Stacking FBPINN results for the single fidelity level 0 and the first two stacking FBPINN levels. \textbf{Right:} Line plots of the results from the stacking FBPINN at $t=0.25$ (top) and $t=0.75$ (bottom).}
    \label{fig:AC}
\end{figure}

We solve the Allen-Cahn equation by dividing the time domain into subdomains, as presented in~\cref{section:dd_time}. The corresponding results for the stacking FBPINN are shown in~\Cref{fig:AC}. The relative $\ell_2$ error for applying two levels of the stacking FBPINN is $5.9{\cdot}10^{-3}$. Previously reported values for the relative error in literature include $1.68{\cdot}10^{-2}$ for the backward compatible PINN \cite{mattey2022novel} and $2.33{\cdot}10^{-2}$ for PINNs with adaptive resampling~\cite{wight2020solving}.

\section{Extension to DeepONets}\label{section:results-deeponets}

\begin{figure}[h!]
    \centering
    \begin{subfigure}[t]{0.65\textwidth}
    	\centering
        \includegraphics[width=\textwidth]{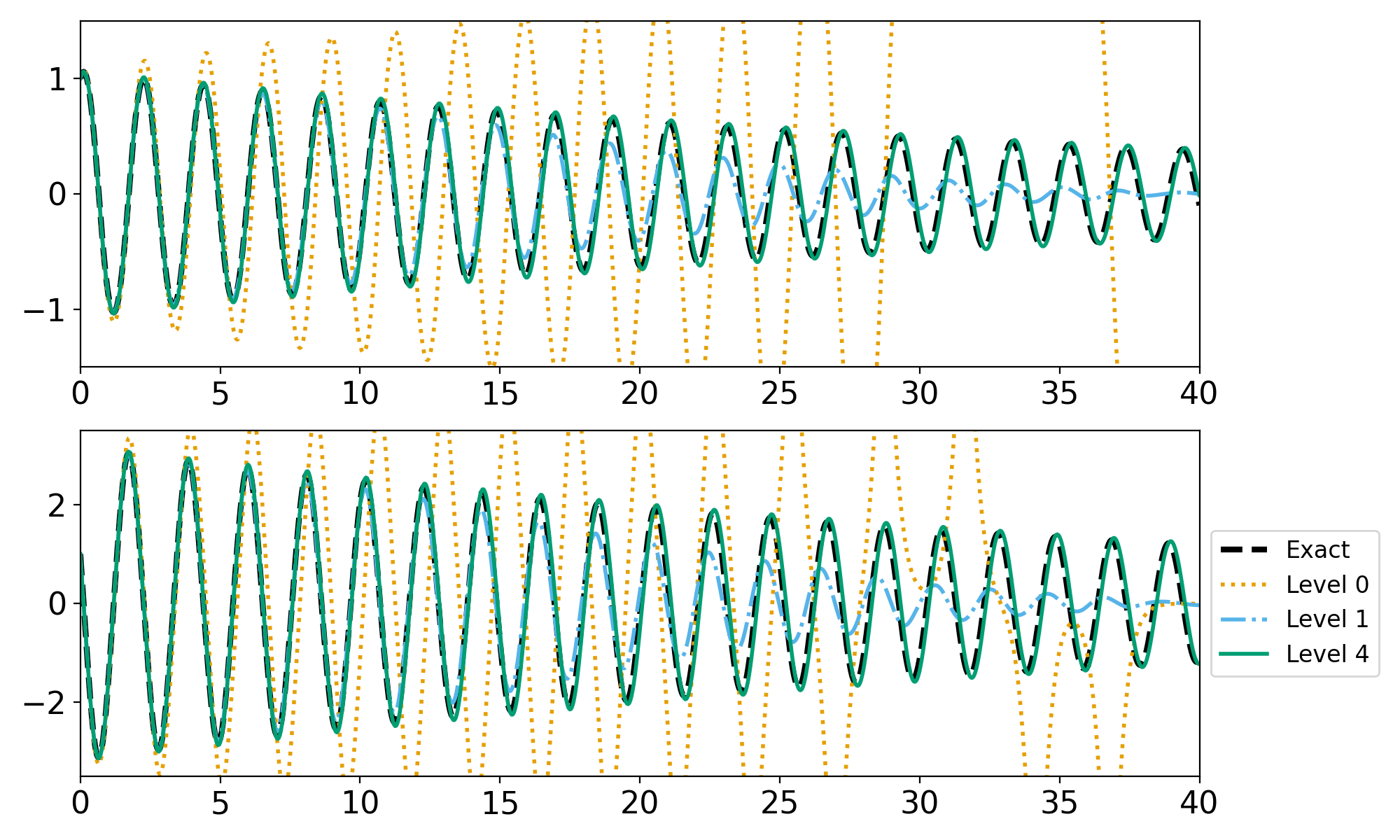}
        \caption{$s_1(0)=1, s_2(0)=1$}
    \end{subfigure}\\
    \begin{subfigure}[t]{0.65\textwidth}    
    	\centering
        \includegraphics[width=\textwidth]{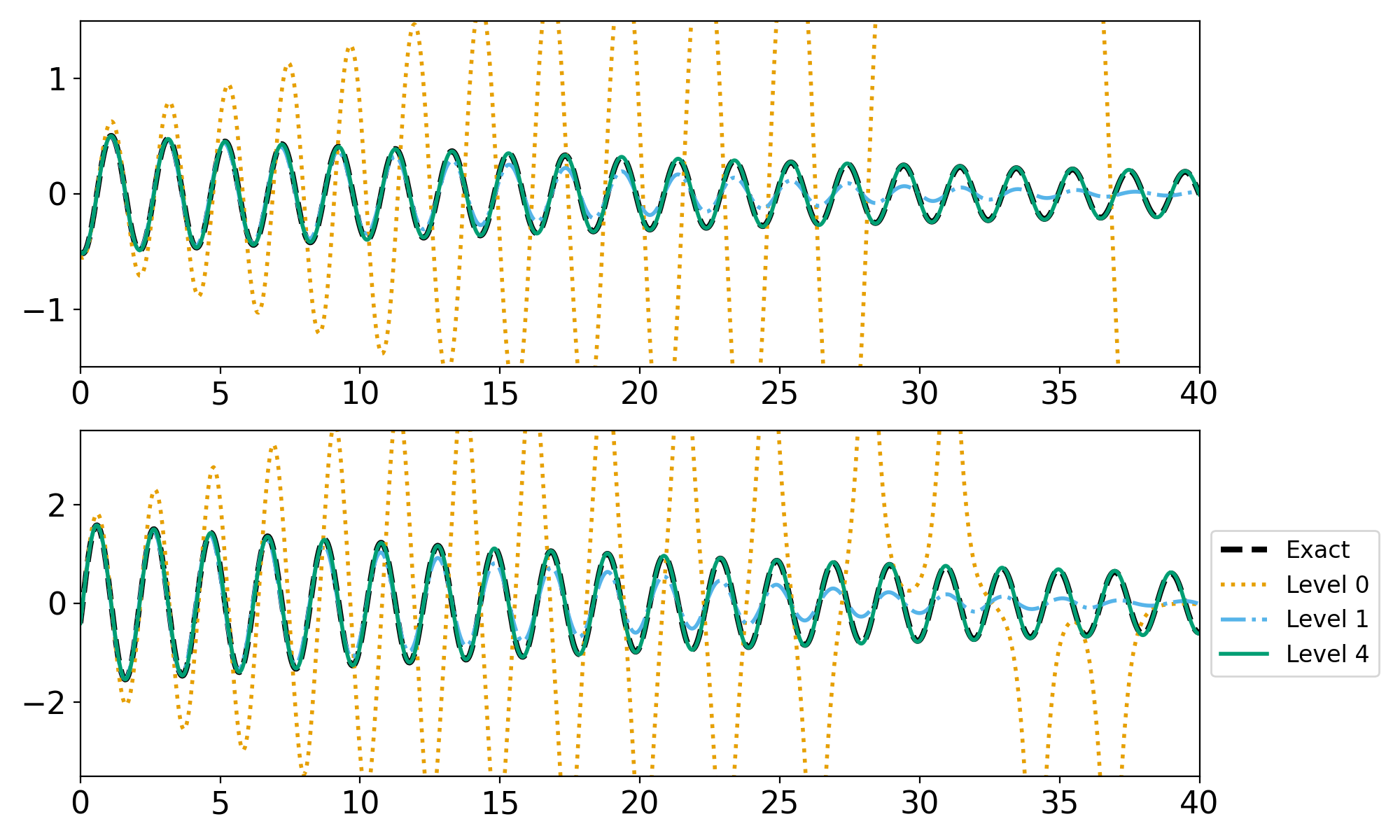}
        \caption{$s_1(0)=-0.5, s_2(0)=-0.4$}
    \end{subfigure}
    \caption{Stacking FB-DON results for the pendulum system at two different sets of initial conditions.}
    \label{fig:FBDeepONet}
\end{figure}

The method presented in~\Cref{section:method} can be extended seamlessly to multifidelity stacking DeepONets from \cite{howard_multifidelity_2023, howard_stacked_2023}; we denote the resulting method as finite-basis DeepONets (FB-DONs). For the sake of brevity, we refer to~\cite{lu_deeponet_2021,howard_multifidelity_2023,howard_stacked_2023} for details on the DeepONet approach. As an example, we present results for the pendulum problem in~\Cref{sec:pendulum} and train a model mapping given initial conditions $\left(s_1(0), s_2(0)\right)$ to the corresponding solution $\left(s_1(t), s_2(t)\right)$ on the whole time interval $[0,20]$. This is referred to as operator learning since we learn a mapping between the initial conditions and the solution space instead of a single solution. One each level $l$, $l > 0$, we train $2^l$ DeepONets with partition of unity functions as defined in~\Cref{eq:PoU}. As training data we employ $50\,000$ randomly chosen pairs $(s_1(0), s_2(0)) \in [-2, 2] \times [-1.2, 1.2]$, and the loss is given by~\Cref{eq:loss} and the differential equations in~\Cref{eq:pendulum_1,eq:pendulum_2}. After training, the resulting FB-DON model is then able to predict the solution for any initial condition in the training range, as shown in~\Cref{fig:FBDeepONet}.

\section{Discussion}\label{section:discussion}
In this paper, we have introduced the stacking FBPINN and FB-DON approaches. For the considered time-dependent problems, stacking FBPINNs yielded more accurate results than stacking PINNs alone and, in some cases, it additionally required fewer total trainable parameters. This indicates that a domain decomposition in time can greatly improve the performance of stacking PINNs. In contrast to prior work on stacking PINNs and DeepONets, stacking FBPINNs and FB-DONs use a sum of subdomain networks weighted by the partition of unity functions on the corresponding level. In contrast to multilevel FBPINNS in~\cite{dolean_multilevel_2023}, in which the subdomain networks are summed across all levels and trained simultaneously, the architecture and training of stacking FBPINNs and FB-DONs is sequential with respect to the levels; the idea is similar to multiplicative coupling as discussed in~\cite{dolean_finite_2022} but implemented differently using the stacking approach. This difference allows for stacking FBPINNs and FB-DONs to consider different equations on different levels, akin to simulated annealing, as considered in~\cite{howard_stacked_2023}, or to consider different physical models at different length scales. We leave this for future work. The extension to stacking FB-DONs allows for use of physics-informed FB-DONs as surrogate models in place of traditional numerical solvers. The computation of a solution using a trained stacking FB-DONs is very efficient: it requires only one forward pass of the networks, and therefore, the computational time compared with classical numerical solvers can be greatly reduced. One advantage of the FB-DONs approach is that it can be used in conjunction with existing methods for increasing accuracy of physics-informed DeepONets, including long-time integration \cite{wang2023long} and adaptive weighting schemes \cite{wang2022improved, howard2024conjugate, qadeer2023efficient}.

\section{Acknowledgments}
    This project was completed with support from the U.S. Department of Energy, Advanced Scientific Computing Research program, under the Scalable, Efficient and Accelerated Causal Reasoning Operators, Graphs and Spikes for Earth and Embedded Systems (SEA-CROGS) project (Project No. 80278). Pacific Northwest National Laboratory (PNNL) is a multi-program national laboratory operated for the U.S. Department of Energy (DOE) by Battelle Memorial Institute under Contract No. DE-AC05-76RL01830. The computational work was performed using PNNL Institutional Computing.

\section{Training parameters}
\begin{table}[h]
    \centering
    \begin{tabular}{ l | c | c |c } 
        & \Cref{sec:pendulum} & \Cref{sec:Multiscale} & \Cref{sec:AC}   \\ \hline
Level 0 learning rate \& decay rate &   $5 {\cdot} 10^{-3}, 0.99$ & $10^{-3}, 0.99$ & $10^{-4}, 0.99$ \\
Level 0 network width & $100$  &  $32$ & $100$ \\
Level 0 network layers & $3$ & $3$ & $6$ \\
Level 0  iterations &   $200\,000$   & $200\,000$   &   $200\,000$   \\        
Nonlinear network width & $32$ & $16$ & $200$ \\
Nonlinear network layers & $3$ & $4$ & $4$\\
Linear network size &   $[2,  4, 2]$ & $[1,  5, 1]$ &  $[1, 5, 1]$ \\
MF learning rate \& decay rate  &    $5 {\cdot} 10^{-3}, 0.99$ & 
$5 {\cdot} 10^{-3}, 0.95$ & $5 {\cdot} 10^{-3}, 0.95$ \\
BC batch size &    $1$      & $1$ & $128$ \\
Residual batch size &   $400$ & $400$ & $1024$\\
Iterations &   $200\,000$      & $300\,000$   &       $300\,000$ \\
$\lambda_{r}$, $\lambda_{bc}$, $\lambda_{\alpha}$&   $1.0, 1.0, 1.0$ & $10.0, 1.0, 1.0$ &  $10.0, 1.0, 10^{-5}$ \\  
Level 0 activation function & swish & swish & tanh \\
MF activation function & swish &swish &swish \\
    \end{tabular}
    \caption{Training parameters for the FBPINN results in the paper. The learning rate is set using the \texttt{exponential\_decay} function in Jax \cite{jax2018github} with the given learning rate and decay rate and 2000 decay steps. The training parameters used for the stacking PINN results are given in~\cite{howard_stacked_2023}.}
    \label{tab:parameters_PINNs}
\end{table}

\begin{table}[h]
    \centering
    \begin{tabular}{ l | c } 
        & \Cref{section:results-deeponets}\\ \hline
Level 0 learning rate \& decay rate &   $5 {\cdot} 10^{-3}, 0.9$ \\
Level 0 branch and trunk width & $100$  \\
Level 0 branch and trunk layers & $5$\\
Level 0  iterations &   $100\,000$  \\        
Nonlinear branch and trunk width & $100$  \\
Nonlinear branch and trunk layers & $3$ \\
Linear branch and trunk width & $10$  \\
Linear branch and trunk layers & $1$ \\
MF learning rate \& decay rate  & $5 {\cdot} 10^{-3}, 0.9$    \\
BC batch size & $1\,000$   \\
Residual batch size &  $10\,000$\\
Iterations &   $200\,000$      \\
$\lambda_{r}$, $\lambda_{bc}$, $\lambda_{\alpha}$&  $1.0, 1.0, 1.0$  \\  
Level 0 activation function & sin \\
MF activation function & sin \\
\end{tabular}
    \caption{Training parameters for the FBDeepONet results in the paper. The learning rate is set using the \texttt{exponential\_decay} function in Jax \cite{jax2018github} with the given learning rate and decay rate and 2000 decay steps.  }
    \label{tab:parameters_DONs}
\end{table}

\bibliographystyle{plain}
\bibliography{main}

\end{document}